\documentclass[twocolumn]{article}

\usepackage{arxiv}

\usepackage[utf8]{inputenc} 
\usepackage[T1]{fontenc}    
\usepackage{textcomp}
\usepackage{hyperref}       
\usepackage{url}            
\usepackage{booktabs}       
\usepackage{amsfonts}       
\usepackage{amsmath}
\usepackage{nicefrac}       
\usepackage{microtype}      
\usepackage{cleveref}       
\usepackage{graphicx}
\usepackage[numbers,sort&compress]{natbib} 
\usepackage{doi}
\usepackage{capt-of}

\title{Infectious behaviour: Simulating the effects of communication and social influence on pathogen transmission in crowds}

\date{}

\usepackage{authblk}

\author[1]{Sophia Johanna Wagner\thanks{Corresponding author: \href{mailto:sophia.wagner@hm.edu}{\texttt{sophia.wagner@hm.edu}}}}
\author[2]{Anne Templeton}
\author[1]{Gerta Köster}
\affil[1]{Munich University of Applied Sciences HM, Department of Computer Science and Mathematics, Lothstraße 64, 80335 München, Germany}
\affil[2]{University of Edinburgh, School of Philosophy, Psychology and Language Sciences, Department of Psychology, 7 George Square, Edinburgh EH8 9JZ, UK}

\renewcommand{\shorttitle}{Infectious behaviour}

\hypersetup{
pdftitle={Infectious behaviour: Simulating the effects of communication and social influence on pathogen transmission in crowds},
pdfauthor={Sophia Johanna Wagner, Anne Templeton, Gerta Köster},
pdfkeywords={crowd dynamics, behavioral dynamics, agent-based modeling, social influence, shared social identity, airborne pathogen exposure, risk communication, dynamical systems},
colorlinks=true,
linkcolor=black,
citecolor=black,
}

\begin{document}
\twocolumn[{
  \begin{@twocolumnfalse}

	\maketitle

	\begin{abstract}
    Public health measures at mass gatherings work only if people follow them, and adherence depends both on how instructions are communicated and what others do. We link these behavioural factors to pathogen transmission in a local-scale, agent-based exposure model combining pedestrian dynamics with airborne transmission.
    Rather than fitting a parametric behavioural model, agents draw their adherence at run time from survey respondents (2,170 attendees of UK sports and music events) who reported the same context:
    we vary stewards' communication effectiveness and adherence of role models, as well as agents' shared social identity with each, assuming that strong shared social identity increases their influence.
    In a ticket-checkpoint queue of 101 agents, effective communication combined with strong shared social identity with stewards reduces the number of highly exposed agents by 82\% for mask wearing. In contrast, strong shared social identity with non-adhering role models gives almost nine times as many highly exposed agents as with adhering role models. Physical distancing was not monotonically beneficial, because adherence changes how agents move: an infectious agent that kept distance moved along the edge of the queue, leading to similar results across conditions.
    Our pipeline transfers to any survey that includes behavioural measures and contextual variables.

	\end{abstract}

  \end{@twocolumnfalse}
  \vspace{1em}
}]

\section{Introduction}

During the COVID-19 pandemic, many superspreading events occurred at crowded events such as sports matches. For example, four days after a professional indoor basketball match in Germany in 2020, 36 out of 69 attendees tested positive~\cite{pauser-2021-life}.
To make these events safer, organisers relied on public health measures such as mask wearing and physical distancing. Whether these measures work, however, depends on whether people follow them. This in turn depends on factors such as how clearly an instruction was given or what other people in the room are doing.

Infection models offer a way to assess the effectiveness of such measures.
Most of them work at the level of whole populations. The classic example is the compartmental SIR model, which divides the population into susceptible, infected, and recovered compartments and moves proportions between them at fixed rates over weeks or months~\cite{friasmartinez-2011-cdyn}. Everyone inside a compartment behaves the same way, so differences between individuals disappear~\cite{friasmartinez-2011-cdyn}.
Researchers have incorporated human behaviour in several ways: letting the transmission rate respond to the number of nearby infections~\cite{eksin-2019-cdyn}, weighing the costs and benefits of a measure~\cite{saadroy-2023-cdyn}, or adding compartments for vaccination~\cite{seibel-2025-cdyn} and for protective habits that spread through social contacts~\cite{ryan-2024-cdyn}.

A recurring challenge is to parameterise behavioural models~\cite{ryan-2024-cdyn}. Many studies use survey data and fit statistical models to the responses, which then determine agent behaviour.
De Gaetano et al.~\cite{degaetano-2024-cdyn} fitted logistic functions to the survey data to determine transitions between compliant and non-compliant compartments and derived state-specific contact matrices that influence disease transmission.
Rodriguez-Cartes et al.~\cite{rodriguez-cartes-2024-cdyn} also used survey data, but in an agent-based network model in which nodes represent individuals and edges their connections. They fitted a logistic regression model to predict an agent's decision to wear a mask from factors such as perceived severity and perceived benefits in a network of 10,000 people.
Durham et al.~\cite{durham-2012-cdyn} followed a similar approach for influenza, deriving each agent's decision to vaccinate or avoid crowds from survey responses about risk perception and beliefs.

None of the models above are designed to capture what happens in a single room over a few minutes or hours, where the exact movement of people and the distance between them decide who gets exposed. This is the domain of local-scale infection models. Balkan et al.~\cite{atamer-balkan-2024-cdyn}, for example, combine a pedestrian model with a grid-based model of virus spread to study measures such as mask wearing.
Ghoroghi et al.~\cite{ghoroghi-2022-cdyn} couple a computational fluid dynamics simulation with an agent-based model to test mask wearing and vaccination at different, but fixed, compliance levels. Most local-scale models, however, do not draw on behavioural data at all when deciding how agents adhere.

In our own previous work~\cite{wagner-2026-cdyn}, we investigated varying levels of adherence to mask wearing and physical distancing within a local-scale exposure model. We also sampled adherence values directly from survey responses, in order to introduce heterogeneous behaviour. However, adherence did not depend on the social and contextual factors that shape it.
To address this, we build on a survey by Templeton et al.~\cite{templeton-2021-cdyn} of 2,170 attendees at sports and music events in the United Kingdom. They were asked about factors associated with their adherence to public health measures. The findings highlight two factors: unclear communication of safety measures is associated with lower adherence, and seeing others adhere is associated with higher adherence.

A study on hospital workers also showed for that clear communication was positively associated with higher self-reported adherence to COVID-19 safety measures~\cite{hlubek-2022-life}.
The influence of other people's behaviour and the way instructions are received point to a psychological mechanism: shared social identity.
People with a shared social identity perceive themselves as part of the same group, which incurs a relational transformation of them feeling more positive and trusting of each other, and increased social influence~\cite{neville-2020-life, reicher-2010-life}.
Neville et al.~\cite{neville-2021-life} report that people adhere more to health measures when they identify with the group around them, when they see others around them adhering, and when instructions come from a person who is seen as one of the group. In line with this, adherence to COVID-19 safety measures was shown to increase when hospital workers perceived their leaders as members of the same group~\cite{smith-2024-cdyn}.

However, shared social identity can also reduce risk perception and therefore increase risk-taking~\cite{cruwys-2020-life}. In support of this, Cruwys et al.~\cite{cruwys-2021-life} found that people with a shared social identity perceived a lower risk of contagion than those outside their group.

\begin{center}
    \includegraphics[width=\columnwidth]{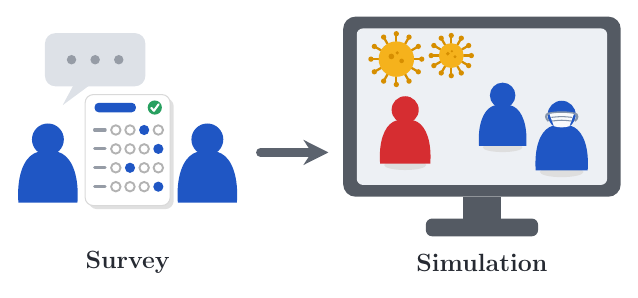}
    \captionof{figure}{Idea of this work: integrating survey data on adherence to public health measures and associated factors into a pathogen transmission model.}
    \label{fig:graphical_abstract}
\end{center}

In this work, we therefore quantify the effect of communication and others' adherence on infection spread, each in combination with shared social identity.
For this, we integrate survey data into an agent-based exposure model, as depicted in Figure~\ref{fig:graphical_abstract}.
We use the local-scale exposure model of Rahn et al.~\cite{rahn-2022-cdyn}, which extends a pedestrian dynamics model with airborne infection transmission, and apply it to a scenario in which 101 agents, including one infectious agent, move towards a ticket checkpoint.
While previous approaches have usually fitted a statistical model to the responses~\cite{ryan-2024-cdyn, rodriguez-cartes-2024-cdyn, durham-2012-cdyn}, thereby assuming a particular relationship between the variables, we sample directly from the survey responses.
This preserves the full spread of responses that a fitted curve would smooth away, and avoids assumptions about the shape of the underlying relationship.
In contrast to our previous work~\cite{wagner-2026-cdyn}, adherence is no longer a fixed property of an agent: it is drawn at run time from those survey respondents who reported the same social or communicative context as the agent perceives.

We aim to answer the following research questions:
\begin{itemize}
 \item How can we integrate survey responses into a local-scale exposure model to quantify the effect of human behaviour on pathogen transmission?
 \item How do communication effectiveness and shared social identity with the person giving instructions influence infection risk?
 \item How do others' adherence and shared social identity with others influence infection risk?
\end{itemize}

Section~\ref{sec:methods} introduces the survey and the simulation framework, Section~\ref{sec:workflow} describes how the two are coupled, and Sections~\ref{sec:results_communication} and~\ref{sec:results_social} report the two simulation studies in which we model the effects of communication and social influence on pathogen transmission.
Finally, we discuss our results in Section~\ref{sec:discussion} and conclude in Section~\ref{sec:conclusion}.

\section{Materials and methods}
\label{sec:methods}

\subsection{Survey data}
\label{sec:survey}

The behavioural input for our model is drawn from a survey by Templeton et al.~\cite{templeton-2021-cdyn}, carried out under the UK Government's 2021 Events Research Programme. The survey gathered responses from 2,170 people attending a range of sports and music events across England, such as the FA Cup Semi-Final and Final, the Carabao Cup Final, the World Snooker Championship, and a music event at Sefton Park. Its aim was to characterise how attendees experienced these events, how they perceived the COVID-19 guidance in place, and which factors shaped their self-reported adherence to public health measures such as physical distancing and mask wearing.
To complement the quantitative responses, the authors also interviewed 37 attendees to examine the reasoning behind their behaviour.
Independent observational data from Gwynne et al.~\cite{gwynne-2024-cdyn} corroborated the self-reported behaviour as part of the Events Research Programme.

From this survey we take the two mechanisms introduced above, communication and adherence of others, each combined with shared social identity. Attendees who were unsure why a measure was needed, how to follow it, or when it applied were more likely not to comply, which points to the role of effective communication. Observing others adhere was associated with higher self-reported adherence, whereas observing non-adherence was associated with lower adherence. The specific survey items used for each mechanism are given in Sections~\ref{sec:comm_survey} and~\ref{sec:social_survey}.
Throughout, we interpret an attendee's self-reported adherence at a real event as the probability that a simulated agent follows a single instruction.

\subsection{Exposure model}
\label{sec:exposure_model}

All simulations in this work are carried out in Vadere~\cite{kleinmeier-2019-cdyn}, an open-source software framework for pedestrian dynamics. It simulates agents moving towards targets while they keep their distance from obstacles and other agents.
Agents move according to the optimal steps model~\cite{seitz-2012-cdyn}: in continuous space, an agent moves in an event-driven fashion to the position that maximises a utility function, attracted by its target and repelled by obstacles and by other agents. The model has been validated against controlled experiments and real-world observations~\cite{seitz-2012-cdyn, seitz-2016c-cdyn, sivers-2016b-cdyn}.

Rahn et al.~\cite{rahn-2022-cdyn} added an exposure model on top of this framework to simulate airborne transmission of infectious diseases in small-scale scenarios, such as restaurants or waiting queues. Every agent is either infectious or susceptible. With each exhalation, an infectious agent releases a stationary aerosol cloud carrying an initial pathogen load. That load decays exponentially over time, while the radius of the cloud grows. Where clouds overlap, their loads add up. A susceptible agent positioned inside an aerosol cloud inhales pathogens at a certain rate, so that we can measure its accumulated exposure.

To interpret these exposure values, Rahn et al.~\cite{rahn-2022-cdyn} defined a high-risk benchmark. According to the Robert Koch Institute, a person is at high risk of infection after spending at least ten minutes within 1.5 metres of an infectious individual~\cite{rki-2022b-life}. Simulating this reference situation yields an exposure value that we use throughout as the threshold above which an agent counts as highly exposed. The model has been validated against documented superspreading events, namely a restaurant outbreak~\cite{lu-2020-life} and a choir rehearsal~\cite{miller-2020-life, hamner-2020-life}.

Physical distancing behaviour was introduced in Vadere by Mayr et al.~\cite{mayr-2021-cdyn}. The parameter \textit{socialDistance} specifies the minimum distance an agent aims to keep from others, which the locomotion model realises as a repulsive potential: the force pushing an agent away from its neighbours grows once they come closer than this preferred distance. In congested situations, the intended distance may be ignored if agents get stuck.

In our previous work~\cite{wagner-2026-cdyn}, we extended the exposure model with mask wearing. Masks are represented by two filtration efficiencies, one applied at exhalation and one at inhalation. A mask worn by an infectious agent reduces the pathogen load of the emitted cloud by the exhalation filtration efficiency, while a mask worn by a susceptible agent reduces the amount of pathogen absorbed by the inhalation filtration efficiency.

A further Vadere component we rely on is the psychology layer of Kleinmeier et al.~\cite{kleinmeier-2020-cdyn}, which represents the cognitive processes underlying human behaviour. First, agents perceive stimuli such as an instruction or a perceived threat in the perception layer. The stimuli are then interpreted and possibly augmented with further information in the cognition layer, before agents carry out the resulting behaviour in the behavioural layer. Depending on the scenario, this may mean waiting, selecting a different target, or, as in the present study, wearing a mask or keeping physical distance. The layer thus models the sequence of perceiving, processing, and acting that underlies human decision-making.
It has been used successfully to reproduce empirically observed human behaviour in different contexts~\cite{kleinmeier-2021-cdyn}.

The simulation workflow combines Python-based data processing with the Java-based crowd simulation framework Vadere~\cite{kleinmeier-2019-cdyn}. A Python script samples adherence values from the survey data and provides the resulting agent parameters as input to Vadere. The package \textit{flowcontrol}~\cite{flowcontrol-2026} sends simulation and survey data between Vadere and the Python script, while the \textit{suq-control} package~\cite{suqcontrol-2026} manages and runs multiple simulations in parallel. 

\section{Simulation workflow}
\label{sec:workflow}

In our model, agents receive an instruction to wear a mask or to keep physical distance. Each simulation studies one of the two measures at a time. Whether an agent follows the instruction depends on its adherence, a probability we draw from survey data. Adherence varies with the context in which the instruction is given, for example very effective communication together with strong shared social identity with stewards, which we refer to as a condition.

\begin{figure}[ht]
\centering
\includegraphics[width=\columnwidth]{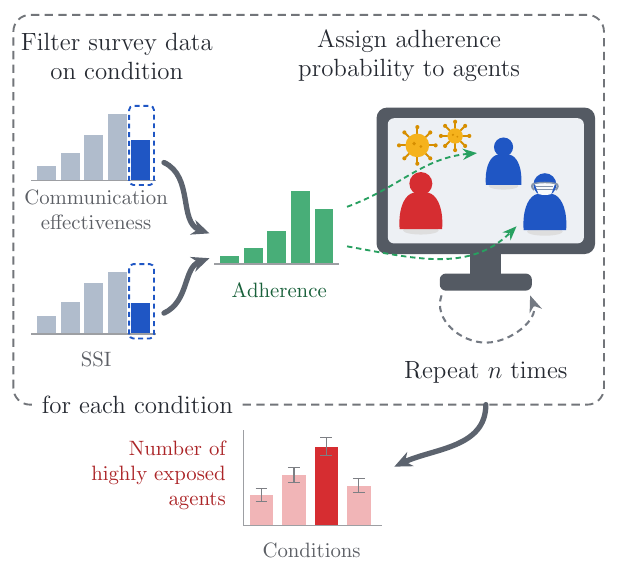}
\caption{Workflow of our model. For a given condition, such as high communication effectiveness combined with strong shared social identity (SSI) with stewards, we filter the survey data to obtain a distribution of adherence probabilities. During the simulation, each agent is assigned an adherence probability drawn from this distribution. Because of the randomness inherent in the model, we repeat this process $n$ times per condition, and finally compare the number of highly exposed agents across conditions.}
\label{fig:workflow}
\end{figure}

Figure~\ref{fig:workflow} summarises the workflow. For a given condition, we select the matching subset of survey respondents, draw adherence values from that subset, and assign them to the agents in the simulation. We repeat each condition many times and compare the resulting number of highly exposed agents between conditions.

\subsection{Preparing the survey data}
\label{sec:workflow_survey_data}

We discretise each input variable into five equal-width intervals, since the variables are measured on a five-point Likert scale. A variable can be the mean of several Likert responses and can therefore also take intermediate values, so we use equal-width interval binning~\cite{dougherty-1995-math}. Crossing two variables yields a $5 \times 5$ grid, and each condition corresponds to one cell of this grid.

Selecting the cell for a condition leaves a subset of respondents with an associated distribution of adherence values. Because these values are also reported on a five-point Likert scale, we convert them to adherence probabilities between 0 and 1, as in our previous work~\cite{wagner-2026-cdyn}. Treating Likert responses as probabilities for computational purposes follows earlier crowd simulation studies that integrate survey data~\cite{mayr-2023-cdyn}. This assumes equidistant psychological intervals between scale points.

The $5 \times 5$ grid produces cells of very unequal size, and some cells contain too few respondents to provide a sufficiently well-supported empirical distribution of adherence values. This is visible where a sparse cell's mean adherence deviates sharply from that of its neighbours. We therefore require at least 30 respondents per cell. We use the minimum denominator sample size specified in the National Center for Health Statistics (NCHS) data-presentation standards for proportions as a pragmatic threshold for the minimum size of a subgroup~\cite{parker-2017-math}. Cells below this size are enriched using a cell-based donor scheme, similar to Andridge and Little~\cite{andridge-2010-math}: a sparse cell absorbs the observations of its four ordinally adjacent neighbours, then its four diagonal neighbours, and successive outward rings as needed, until the pooled count reaches 30 (see Supplementary material 1 for the grid before and after enrichment).

\subsection{Modelling agent behaviour using survey responses}

\begin{figure}[ht]
\centering
\includegraphics[width=\columnwidth]{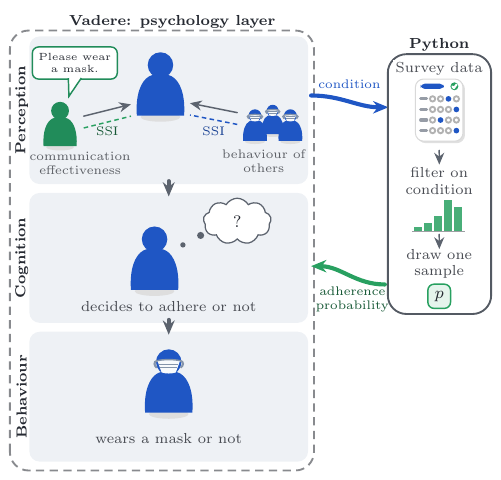}
\caption{Decision process of a single agent for mask wearing. In the perception layer the agent receives an instruction together with the factors under study (communication effectiveness and shared social identity (SSI) with the stewards, or the adherence of others and the shared social identity with them), which define its condition. Through the \textit{flowcontrol} package this condition is sent to a Python script, which samples an adherence probability from the matching survey cell and returns it. In the cognition layer the agent turns this probability into a binary decision, and in the behavioural layer it implements the decision.}
\label{fig:psychology_layer}
\end{figure}

An agent's adherence probability is not fixed before a run. Instead, each agent obtains it at the moment it receives an instruction, so its behaviour reflects the condition it finds itself in. This process runs through the psychology layer described in Section~\ref{sec:exposure_model}. Figure~\ref{fig:psychology_layer} shows the process for a single agent receiving an instruction to wear a mask.

In the perception layer, an agent receives a stimulus: an instruction to wear a mask
or to keep physical distance, as implemented in our previous work~\cite{wagner-2026-cdyn}. We enrich this stimulus with the factors under study. An instruction can carry a communication effectiveness, and each agent can carry a shared social identity with the instruction giver or with the surrounding crowd, both newly added attributes. Depending on the study, the agent also perceives whether the agents around it adhere. Together, these perceived factors define the agent's condition, for example an instruction perceived as very effective from a steward the agent has a strong shared social identity with.

We use \textit{flowcontrol}~\cite{flowcontrol-2026} to send this condition to a Python script, where it selects the corresponding cell of the survey grid. We draw one adherence probability at random from that cell's distribution and return it to the agent, which is assigned this value. In the cognition layer, the agent then turns the probability into a binary decision: it draws a uniform random number and adheres if the number falls below its adherence probability. Finally, in the behavioural layer, the agent implements its decision, wearing a mask or keeping distance, or doing neither. Because the probability is assigned at run time based on what the agent perceives, agents in the same run can behave differently, and the same agent would behave differently under a different condition.

The adherence of the infectious agent is not sampled in this way. It is set by the condition under study, so that its effect on the number of highly exposed agents is controlled rather than random.

\subsection{Running multiple simulations}

The simulation is stochastic in three ways: adherence values are sampled at random, the decisions taken from them are random, and agent movement is stochastic, so the position of the infectious agent in the queue varies between runs. We therefore run each condition $n$ times with different random seeds and report the mean and standard deviation of the number of highly exposed agents, as defined in Section~\ref{sec:exposure_model}.

Following Rahn~\cite{rahn-2024c-cdyn}, we determine $n$ by visually comparing distributions of simulation outputs rather than by applying a formal stopping rule, which might over- or underestimate $n$ for skewed outputs. We run each condition with $n = 10^3$ and compare the distribution against that of a subset of $n = 10^2$ runs (see Supplementary material 2). The deviation is small in all conditions, so we consider $n = 10^3$ repetitions representative.

\section{Modelling the effect of communication on pathogen transmission}
\label{sec:results_communication}

When a steward in a stadium entrance queue asks the crowd to wear masks or keep physical distance, individual compliance may depend on how clearly the instruction is understood, how important the measure is perceived to be, and whether the steward is seen as part of one's own group. 
In this section we vary two factors, how effective the communication is perceived and how much shared social identity people feel with the stewards, and quantify how they affect the spread of infection.

\subsection{Simulation setup}

\begin{figure*}[ht!]
\centering
\includegraphics[width=\textwidth]{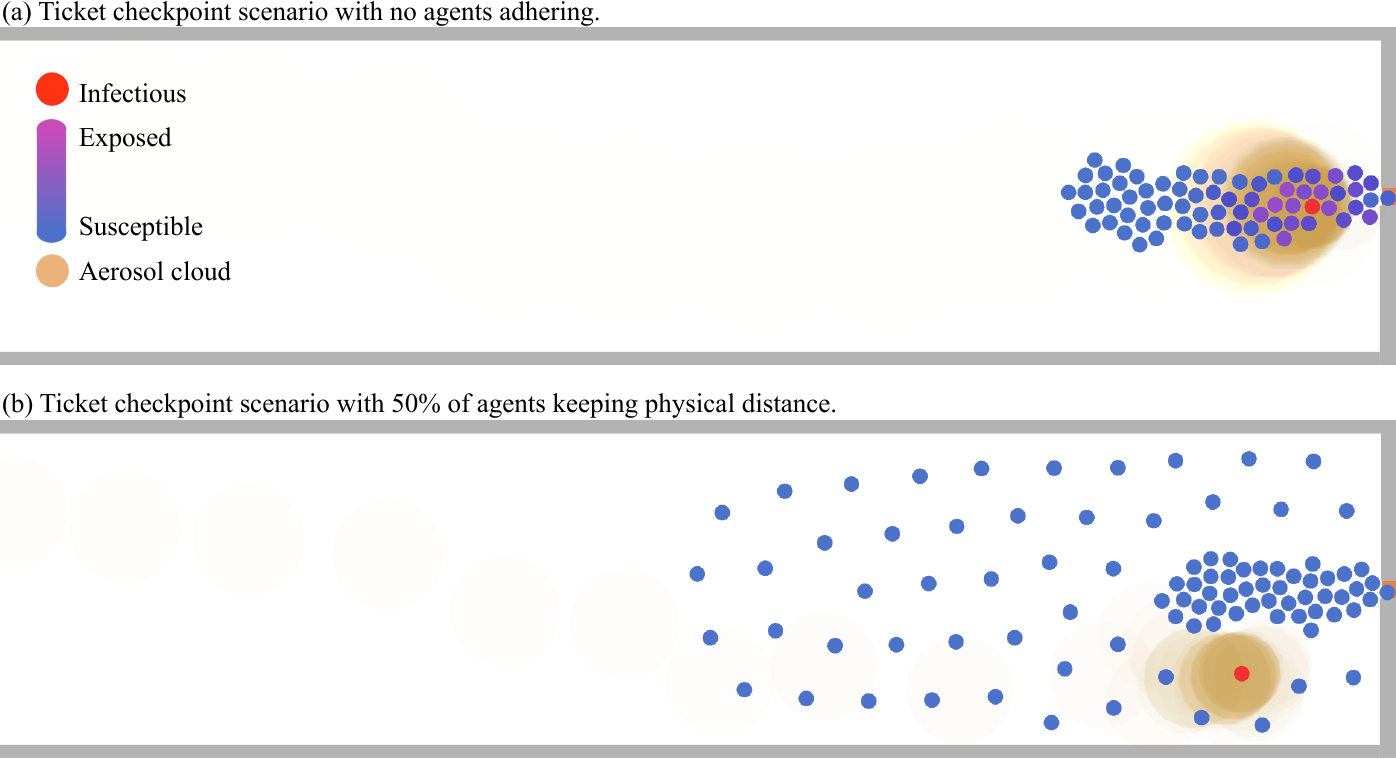}
\caption{Ticket-checkpoint scenario in Vadere. Agents move towards the ticket checkpoint on the right. (a) No agents keep physical distance. (b) Half of the agents keep physical distance.}
\label{fig:vadere_scenario}
\end{figure*}

We use the ticket-checkpoint scenario from our previous work~\cite{wagner-2026-cdyn}. One infectious and 100 susceptible agents move in an unstructured queue towards a checkpoint (Figure~\ref{fig:vadere_scenario}). Each agent is served in 15 seconds, which gives a total simulation time of 25 minutes. 

We consider two public health measures, mask wearing and physical distancing. For physical distancing, agents are instructed to keep a distance of 1.5 metres, following the recommendation of the Robert Koch Institute during the SARS-CoV-2 pandemic~\cite{ibbs-2020-life}. The Events Research Programme behind our survey prescribed 2 metres instead. We tested both and found the difference between them negligible, so we report the 1.5-metre results throughout. When agents keep their distance, the queue spreads out and fewer agents stand within the trail of an aerosol cloud.

For mask wearing, we consider cases in which the infectious agent does not wear a mask. In our previous work~\cite{wagner-2026-cdyn} the number of highly exposed agents dropped to zero whenever the infectious agent wore a mask, so a mask-wearing infectious agent leaves nothing to compare between conditions.
We set the inhalation filtration efficiency to 89\%. This is the conservative minimum guaranteed for FFP2 respirators by the European standard EN~149:2001+A1:2009~\cite{europeanstandard-2009-life}.
The parameters we use for the exposure model were validated by Rahn~\cite{rahn-2024c-cdyn} and are listed in Supplementary material 3.

\subsection{Survey data analysis}
\label{sec:comm_survey}

The two survey variables for this study are the communication effectiveness of the stewards and the shared social identity with the stewards. Communication effectiveness comes from a single question, in which attendees rated how effective they found the stewards' communication at the event on a five-point Likert scale, from very ineffective to very effective.
Shared social identity is the mean of three questions, in which attendees rated how similar to the stewards they felt, how much unity they felt with them, and how much togetherness there was between attendees and stewards. The five-point Likert scale ranged from strongly disagree to strongly agree. In this work, we refer to a perceived shared social identity of ``strongly disagree'' as weak and ``strongly agree'' as strong.

\begin{figure}[ht]
\centering
\includegraphics[width=\columnwidth]{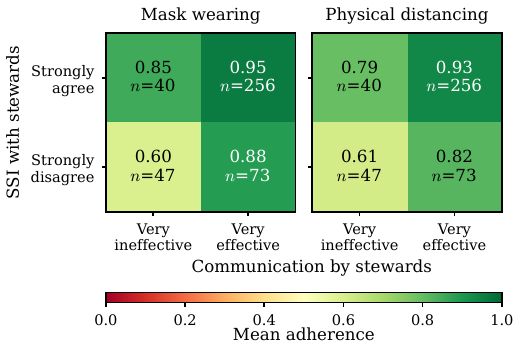}
\caption{Mean adherence and respondent counts for the corner cells of the survey grid, for communication effectiveness of the stewards and shared social identity (SSI) with the stewards. Adherence is shown for mask wearing and physical distancing, after iterative merging of sparse cells to a minimum of 30 respondents (Section~\ref{sec:workflow_survey_data}). Mean adherence is shown for readability, but the simulation samples from the full distribution within each cell.}
\label{fig:comm_survey_data}
\end{figure}

We compare the corner cells of the grid: strong versus weak shared social identity, and very effective versus very ineffective communication. Figure~\ref{fig:comm_survey_data} shows the mean adherence for mask wearing and physical distancing in these cells, together with the number of respondents in each. Adherence is lowest when shared social identity is weak and communication very ineffective, at 0.60 for mask wearing and 0.61 for physical distancing, and highest when shared social identity is strong and communication very effective, at 0.95 and 0.93. The two factors thus point in the same direction: better communication and stronger shared social identity with the stewards both go together with higher adherence, and their combination gives the highest adherence of all. See Supplementary material 1 for the full $5 \times 5$ grid.

\subsection{Simulation results}

\begin{figure}[ht]
\centering
\includegraphics[width=\columnwidth]{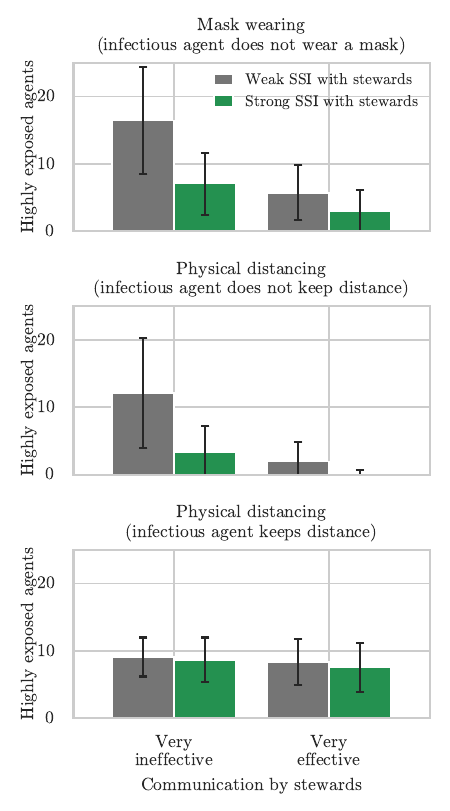}
\caption{Mean and standard deviation (vertical black lines) of the number of highly exposed agents, for mask wearing and for physical distancing with and without the infectious agent keeping distance, across the four combinations of very weak and very strong shared social identity (SSI) with stewards and very ineffective and very effective communication by stewards.}
\label{fig:results_highrisk_communication}
\end{figure}

Figure~\ref{fig:results_highrisk_communication} shows the number of highly exposed agents for each condition for mask wearing, where the infectious agent does not wear a mask, and for physical distancing, with and without the infectious agent keeping distance.
For mask wearing, the number of highly exposed agents is highest when communication is ineffective and shared social identity with the stewards weak, at 16.45 agents, and lowest when both are strong, at 2.98. That is a reduction of 82\%. Each factor alone gives a smaller but clear improvement, which matches the survey data: what raises adherence in the survey lowers the number of highly exposed agents in the simulation.

Physical distancing behaves differently, and the outcome depends on the infectious agent. When the infectious agent does not keep distance, the pattern resembles mask wearing but is sharper: the number of highly exposed agents falls from 12.08 to 0.05 when both communication and shared social identity are strong. The near-zero highly exposed agents can be explained by a pattern we already observed in our previous work~\cite{wagner-2026-cdyn}: when most susceptible agents adhere but the infectious agent does not, it bypasses the queue, which reduces the time it spends near other agents.
When the infectious agent does keep distance, the effect of both factors almost disappears. The number of highly exposed agents stays between 7.50 and 9.06 across all four conditions.
Figure~\ref{fig:vadere_scenario}(b) shows why: when the infectious agent does adhere, it moves out of the inner part of the queue, where the non-adhering agents stand very close to each other. It is therefore never close to many agents at once, but still spends enough time in the queue to leave a trail of aerosol clouds.

The standard deviations are large throughout, in several conditions close to the mean. This is expected: adherence is drawn at random for every agent, and the position of the infectious agent in the queue differs between runs, so single runs vary widely even under the same condition.
See Supplementary material 2 for the distributions of highly exposed agents, and Supplementary material 4 for the means and standard deviations.

\section{Modelling the effect of social dynamics on pathogen transmission}
\label{sec:results_social}

Adherence may also be shaped by the behaviour of others in the crowd. Spectators who see fellow fans of their own team wearing masks may do the same, whereas non-adherence by the same in-group may encourage them to ignore the instruction. The behaviour of fans of the opposing team may influence them less. Here we vary whether another group at the front of the queue adheres and how much shared social identity people feel with that group.

\subsection{Simulation setup}

\begin{figure*}[ht!]
\centering
\includegraphics[width=\textwidth]{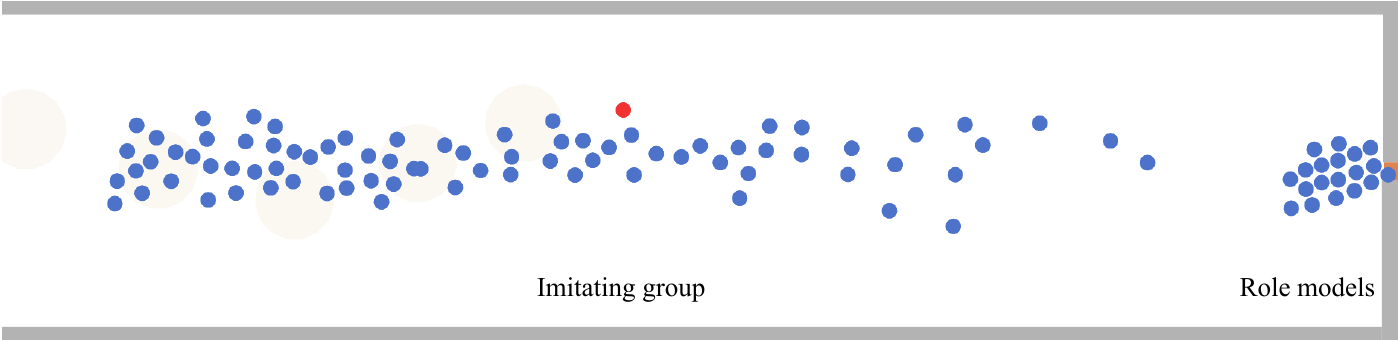}
\caption{Two-group ticket-checkpoint scenario used to study social influence. Twenty role models enter first, followed 30 seconds later by an imitating group of 80 agents and one infectious agent. Each imitating agent decides whether to adhere on the basis of the adherence of the role models ahead and the shared social identity it feels with them.}
\label{fig:social_scenario}
\end{figure*}

We keep the ticket-checkpoint scenario, the model parameters, and the same two measures, mask wearing and physical distancing, as in Section~\ref{sec:results_communication}, but split the crowd into two groups (Figure~\ref{fig:social_scenario}).
First, 20 role models enter and move towards the checkpoint. Thirty seconds later, an imitating group of 80 agents follows, together with one infectious agent. The role models do not draw adherence from the survey: we set them either to adhere or not to adhere, so that they define the social context. The imitating agents decide whether to adhere on the basis of two factors: whether the role models ahead of them adhere, and how much shared social identity they feel with them. 

We do not update adherence during a run. An agent draws its adherence once, from the survey cell matching its condition, and keeps it. Updating it dynamically as agents observe one another would add further assumptions about how and how fast people change their minds, and would tend to push every run towards full adherence, because adherence in the survey data is generally high.

Our quantity of interest here is the number of highly exposed agents in the imitating group only, because these are the agents whose adherence responds to the role models. The results in this study are therefore not directly comparable to those of the communication study, where all susceptible agents are counted.

\subsection{Survey data analysis}
\label{sec:social_survey}

The two survey variables for this study are the adherence of others and the shared social identity with the crowd. The adherence of the role models corresponds to the survey item on adherence of others, specifically how strongly attendees felt that other attendees adhered. There were separate questions for mask wearing and physical distancing, and we use the one matching the measure under study.
The shared social identity with the role models corresponds to the shared social identity with the crowd. This is the mean of three questions on how much attendees felt that everyone in the crowd was part of the same group, felt united, and felt a sense of commonality with one another. All items were rated on a five-point Likert scale from strongly disagree to strongly agree.

\begin{figure}[ht]
\centering
\includegraphics[width=\columnwidth]{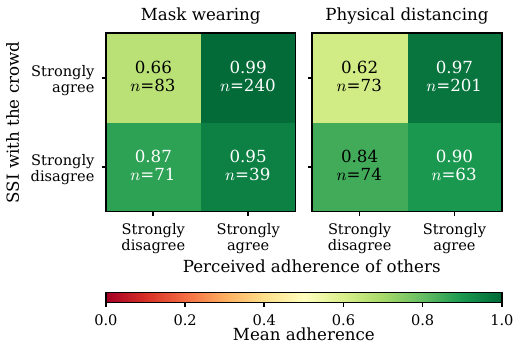}
\caption{Mean adherence and respondent counts for the corner cells of the survey grid, for the adherence of others and shared social identity (SSI) with the crowd. Adherence is shown for mask wearing and physical distancing, after iterative merging of sparse cells to a minimum of 30 respondents. Mean adherence is shown for readability, but the simulation samples from the full distribution within each cell.}
\label{fig:social_survey_data}
\end{figure}

Figure~\ref{fig:social_survey_data} shows the mean adherence for both measures in the four corner cells, with the number of respondents in each. Adherence is highest when the role models adhere and shared social identity is strong, at 0.99 for mask wearing and 0.97 for physical distancing, and lowest when the role models do not adhere and shared social identity is strong, at 0.66 and 0.62.
Shared social identity therefore does not act in one direction: when others adhere, stronger shared social identity raises adherence further, from 0.95 to 0.99 for mask wearing and from 0.90 to 0.97 for physical distancing. When others do not adhere, stronger shared social identity lowers it, from 0.87 to 0.66 and from 0.84 to 0.62. See Supplementary material 1 for the full $5\times 5$ grid.

\subsection{Simulation results}

\begin{figure}[ht]
\centering
\includegraphics[width=\columnwidth]{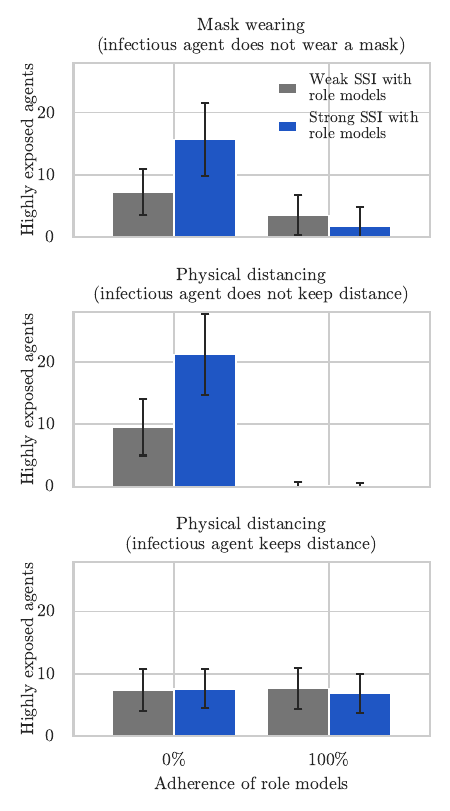}
\caption{Mean and standard deviation (vertical black lines) of the number of highly exposed agents in the imitating group, for mask wearing and for physical distancing with and without the infectious agent keeping distance, across the four combinations of adhering and non-adhering role models with very weak and very strong shared social identity (SSI) with the role models.}
\label{fig:results_highrisk_othersadhere}
\end{figure}

Figure~\ref{fig:results_highrisk_othersadhere} shows the number of highly exposed agents in the imitating group for each condition. The picture differs from the communication study in one important way: shared social identity does not always help.

For mask wearing, when the role models adhere, stronger shared social identity lowers the number of highly exposed agents slightly, from 3.49 to 1.75 agents. When they do not adhere, stronger shared social identity raises it, from 7.25 to 15.65 agents. The worst condition is therefore not the one with the weakest social ties, but the one in which the imitating agents identify strongly with a group that ignores the instruction.

Physical distancing shows the same pattern, and again the infectious agent matters. When the infectious agent does not keep distance, the number of highly exposed agents is near zero if the role models adhere, and rises to 9.54 and 21.23 agents if they do not, for weak and strong shared social identity respectively. As in the communication study, an infectious agent that ignores the instruction while the crowd adheres bypasses the queue and spends less time close to others~\cite{wagner-2026-cdyn}. If the crowd does not keep distance either, that effect disappears, and strong shared social identity with non-adhering role models more than doubles the number of highly exposed agents.
When the infectious agent keeps distance, the differences again almost vanish, with between 6.81 and 7.64 highly exposed agents across all four conditions, for the same reason as in Section~\ref{sec:results_communication}: the adhering infectious agent moves along the edge of the queue rather than through its centre.

As above, the standard deviations are large because of the randomness in the model.
See Supplementary material 2 for the distributions of highly exposed agents, and Supplementary material 4 for the means and standard deviations.

\section{Discussion}
\label{sec:discussion}

We integrated survey data on adherence to public health measures into a local-scale exposure model. Agents in a ticket-checkpoint queue receive an instruction to wear a mask or to keep physical distance. At that moment, each agent draws an adherence probability from the survey respondents who reported the same context. We then compared the number of highly exposed agents across different levels of communication effectiveness and shared social identity with the stewards, and of role model adherence and shared social identity with the role models.
In our approach, the behavioural input is an empirical distribution rather than a fitted function. Sampling within a context cell keeps the full spread of the responses and ties every behaviour to a group of real respondents. In addition, the model translates survey answers into a physical outcome. It turns a Likert score into an agent's decision, and that decision into an inhaled pathogen load, which also depends on what the agents around it do.

In our first simulation study, effective communication lowered the number of highly exposed agents, and shared social identity with the stewards lowered it further. When both were strong, the number of highly exposed agents for mask wearing fell by 82\% compared with very ineffective communication and weak shared social identity with the stewards. This matches the survey authors' recommendation that guidance is best given by people the crowd sees as part of its own group~\cite{templeton-2021-cdyn}.
In related work, Templeton et al.~\cite{templeton-2026-cdyn} found that first responders felt collective language led to the public feeling part of the same group with them and therefore they were more likely to follow the first responders' instructions. We expect a similar effect for public health measures.

In the second simulation study, shared social identity did not always help. When the role models adhered, strong shared social identity with them lowered the number of highly exposed agents slightly, from 3.49 to 1.75. However, when the role models ignored the instruction, strong shared social identity raised it from 7.25 to 15.65 agents. The riskiest condition in our simulations was therefore a crowd with a strong shared social identity with role models who did not adhere.
Thus, in practice, building a sense of shared social identity at an event is not enough on its own. The people the crowd identifies with have to be seen following the guidance.

For mask wearing, higher adherence always means fewer highly exposed agents. For physical distancing it does not, because adherence changes how agents move in the crowd. When the infectious agent kept distance, all four conditions gave between 6 and 10 highly exposed agents, for both the communication and social influence study. In this case, the infectious agent moved along the edge of the queue, never close to many agents at once, but still leaving a trail of aerosol clouds.
When it did not keep distance while the crowd did, the number of highly exposed agents fell to almost zero, because the infectious agent simply bypassed the queue.
This shows that adherence and the number of highly exposed agents are not monotonically related once behaviour affects movement. A model applying a fixed transmission reduction per adhering person cannot reproduce this.

The standard deviations are large, in several conditions close to the mean, despite 1000 runs per condition. We consider this realistic. Adherence is an individual decision, and so is the path an agent takes through the queue. Two events with the same instruction and the same average adherence can therefore still end differently. Because we sample from the full range of survey responses, our model keeps this spread instead of averaging it away.

Several limitations concern the survey data and how we used them.
Self-reported adherence was high in general, partly because the sample included groups recruited for the Events Research Programme, such as healthcare workers, and many participants wanted live events to return~\cite{templeton-2021-cdyn}. A less adherent population would probably increase the overall number of highly exposed agents.
We also assume that the context leads to a certain behaviour. The survey shows which factors occur together, but not which one causes the other. Attendees who adhered may equally well have judged the stewards more favourably afterwards. For the factors we study, the assumed direction is nevertheless the plausible one. An instruction has to be understood before it can be followed, so the communication precedes the decision. The interviews conducted alongside the survey point the same way: attendees who were unsure why a measure was needed or how to follow it gave this as their reason for not adhering~\cite{templeton-2021-cdyn}. 
The survey item on others' adherence asked how strongly attendees felt that others adhered, whereas our role models actually adhere. We assume that the perception matches the actual adherence of others, which plausibly drives the decision.
Finally, we rescaled the five-point answers linearly to the interval between 0 and 1 and used them as adherence probabilities. This assumes that the steps of the scale are equally spaced.

Two further choices follow from the size of the data set. First, we binned each variable into five intervals. Five bins match the five points of the Likert scale, so each bin corresponds to one step a respondent could have chosen. Fewer bins would merge the extreme answers that we compare, and more bins would leave cells too small to sample from. A different number of bins would therefore change the exact values, but not the overall pattern.
Second, we studied communication and social influence separately. Crossing two variables already yields 25 cells, and the respondents are spread unevenly across them, so several cells fall below 30 respondents and have to be enriched from their neighbours. Crossing three variables would yield 125 cells and crossing four would yield 625. Most of them would be empty or consist largely of donated observations, and the conditions would no longer be distinguishable from one another. Quantifying the joint effect of communication and social influence therefore requires a larger survey.

Our model rests on four further assumptions. Aerosol clouds are stationary. Agents in our scenario move slowly and stay in a small area, so this holds here, but it would not hold under fast movement or forced ventilation.
Agents act as individuals, whereas real attendees might arrive as groups of friends or family. Such groups stand close together whatever the instruction, and they may influence each other more strongly than strangers do~\cite{wagner-2026-cdyn}.
Adherence is drawn once per run, whereas real people revise their behaviour in response to situational cues~\cite{alhajri-2024-cdyn}. Updating it would add untested assumptions about the rate of change and, given the high adherence here, would drive nearly every run to full adherence.
In the second study, the imitating agents respond only to the role models, not to each other. The role models enter first and are visible ahead in the queue, so they are the group an arriving attendee can actually observe. Letting agents respond to their immediate neighbours as well would create a cascade whose speed we cannot parameterise in our model. What we report is therefore the effect of the visible reference group alone.

We cannot validate the pipeline end to end, as that would require an outbreak at an event where communication, shared social identity, and adherence were all measured. However, the components are validated individually: locomotion against experiments and field data, exposure against two superspreading events, the psychology layer against observed behaviour, and the survey responses against independent observations of the same events. Their composition remains untested. What we can say is that the simulated behaviour and the resulting numbers of highly exposed agents are plausible in every condition we studied.

The pipeline is not specific to this survey. Our model can draw on any survey responses, given a behavioural measure, at least one contextual variable, and enough respondents per subgroup.
The approach could therefore also be applied to other fields, such as evacuation simulation, where agents react to instructions to leave a building.
Future work could also address more complex scenarios: for example, an attendee moves from a ticket checkpoint to a spectators' area and then to a food court. Each place has its own context, to which the attendee reacts.
Another direction would be to model social groups instead of independent agents. Friends and families arrive together, stand close together, and influence each other in their decisions. This would require survey data on how people behave when they attend an event in groups.

\section{Conclusion}
\label{sec:conclusion}

We presented an approach that turns survey responses directly into agent behaviour in a local-scale exposure model. No parametric behavioural model sits in between. Agents draw their adherence from the respondents who described the same social and communicative context, which keeps the full spread of real answers.

Three messages follow for event organisers. Communication works, but it works best when the crowd identifies with the person speaking, and when the instruction explains why the measure matters and how to follow it. Both together reduced the number of highly exposed agents by 82\% for mask wearing in our scenario.
Social dynamics matter at least as much. Shared social identity amplifies whatever the role model group does. A crowd with a strong shared social identity with role models who ignored the instruction was the worst case in our simulations, with 15.65 highly exposed agents compared with 1.75 when the crowd had a strong shared social identity with role models who did adhere. Creating a sense of unity at an event is therefore not protective on its own. The people the crowd identifies with have to be seen adhering.
Finally, the two measures are not interchangeable. Mask wearing reduced the number of highly exposed agents with increasing adherence. Physical distancing changed how agents moved, and its effect depended on who adhered, above all on the infectious person, which an organiser can not know in advance. Therefore, masks were the more reliable of the two.

\section*{Funding}
The author S.J.W. received funding from the German Research Foundation (DFG) through the project PanVadere (project number 515675334).
The funders had no role in study design, data collection and analysis, decision to publish, or preparation of the manuscript.

\bibliographystyle{unsrtnat}
\bibliography{refs}

\end{document}